\documentclass[12pt,reqno]{amsart}
\usepackage{amsmath,amsfonts,amsthm,amsopn,amssymb}
\usepackage{cite,marginnote}
\usepackage{color,enumitem,graphicx}
\usepackage[colorlinks=true,urlcolor=blue,
citecolor=red,linkcolor=blue,linktocpage,pdfpagelabels,
bookmarksnumbered,bookmarksopen]{hyperref}
\usepackage[english]{babel}

\usepackage[left=2.9cm,right=2.9cm,top=2.8cm,bottom=2.8cm]{geometry}
\usepackage[hyperpageref]{backref}

\numberwithin{equation}{section}

\makeindex
\makeindex
\newtheorem{theorem}{Theorem}[section]
\newtheorem{definition}[theorem]{Definition}
\newtheorem{lemma}[theorem]{Lemma}
\newtheorem{corollary}[theorem]{Corollary}
\newtheorem{proposition}[theorem]{Proposition}
\newtheorem{remark}[theorem]{Remark}

\newcommand{\s}{\section}

\newcommand{\G}{\Gamma}

\newcommand{\R}{\mathbb R}

\newcommand{\C}{\mathbb C}

\newcommand{\lab}{\label}
\newcommand{\bt}{\begin{theorem}}
\newcommand{\et}{\end{theorem}}
\newcommand{\bl}{\begin{lemma}}
\newcommand{\el}{\end{lemma}}
\newcommand{\bd}{\begin{definition}}
\newcommand{\ed}{\end{definition}}
\newcommand{\bc}{\begin{corollary}}
\newcommand{\ec}{\end{corollary}}
\newcommand{\bp}{\begin{proof}}
\newcommand{\ep}{\end{proof}}
\newcommand{\bx}{\begin{example}}
\newcommand{\ex}{\end{example}}
\newcommand{\bi}{\begin{exercise}}
\newcommand{\ei}{\end{exercise}}
\newcommand{\bo}{\begin{proposition}}
\newcommand{\eo}{\end{proposition}}
\newcommand{\br}{\begin{remark}}
\newcommand{\er}{\end{remark}}
\newcommand{\beq}{\begin{equation}}
\newcommand{\eeq}{\end{equation}}
\newcommand{\ba}{\begin{align}}
\newcommand{\ea}{\end{align}}
\newcommand{\bn}{\begin{enumerate}}
\newcommand{\en}{\end{enumerate}}
\newcommand{\bg}{\begin{align*}}
\newcommand{\bcs}{\begin{cases}}
\newcommand{\ecs}{\end{cases}}

\newcommand{\bean}{\begin{eqnarray*}}
\newcommand{\eean}{\end{eqnarray*}}

\def\C{\mathbb{C}}
\def\N{\mathbb{N}}

\def\R{\mathbb{R}}

\def\bd{\mathrm{bd}\,}

\title[Normalized ground state for coupled system]{Normalized ground states  for a coupled Schr\"odinger system: Mass super-critical case}

\author[L.~Jeanjean]{Louis Jeanjean}
\author[J.~J.~Zhang]{Jianjun Zhang}
\author[X.~X.~Zhong]{Xuexiu Zhong}

\address[L.~Jeanjean]{\newline\indent Universit\'e de Franche-Comt\'e
\newline\indent
CNRS, UMR 6623, LmB, F-25000 Besan\c{c}on, France}
\email{\href{mailto:louis.jeanjean@univ-fcomte.fr}{louis.jeanjean@univ-fcomte.fr}}

\address[J.~J.~Zhang]{\newline\indent College of Mathematics and Statistics
\newline\indent
Chongqing Jiaotong University
\newline\indent
Chongqing 400074, PR China}
\email{\href{mailto:zhangjianjun09@tsinghua.org.cn}{zhangjianjun09@tsinghua.org.cn}}

\address[X.~X.~Zhong]{\newline\indent South China Research Center for Applied Mathematics and Interdisciplinary Studies
\newline\indent
South China Normal University
\newline\indent
Guangzhou 510631, PR China}
\email{\href{mailto:zhongxuexiu1989@163.com}{zhongxuexiu1989@163.com}}

\thanks{Xuexiu Zhong NSFC was supported by the NSFC (No.12271184), Guangdong Basic and Applied Basic Research Foundation (2021A1515010034),Guangzhou Basic and Applied Basic Research Foundation(202102020225). Jianjun Zhang was supported by the NSFC (No.12371109)}

\subjclass[2000]{}

\begin{document}

\begin{abstract}
We consider the existence of solutions $(\lambda_1,\lambda_2, u, v)\in \R^2\times (H^1(\R^N))^2$ to systems of coupled Schr\"odinger equations
$$
\begin{cases}
-\Delta u+\lambda_1 u=\mu_1 u^{p-1}+\beta r_1 u^{r_1-1}v^{r_2}\quad &\hbox{in}~\R^N,\\
-\Delta v+\lambda_2 v=\mu_2 v^{q-1}+\beta r_2 u^{r_1}v^{r_2-1}\quad &\hbox{in}~\R^N,\\
0<u,v\in H^1(\R^N), \, 1\leq N\leq 4,&
\end{cases}
$$
satisfying the normalization  $$ \int_{\R^N}u^2 \mathrm{d}x=a \quad \mbox{and} \quad  \int_{\R^N}v^2 \mathrm{d}x=b.$$ Here $\mu_1,\mu_2,\beta>0$ and the prescribed masses $a,b>0$. We focus on the coupled purely mass super-critical case, i.e.,
$$2+\frac{4}{N}<p,q,r_1+r_2<2^*$$
with $2^*$ being the Sobolev critical exponent, defined by $2^*:=+\infty$ for $N=1,2$  and $2^*:=\frac{2N}{N-2}$ for $N=3,4$. We optimize the range of $(a,b,\beta,r_1,r_2)$ for the existence. In particular, for $N=3,4$ with $r_1,r_2\in (1,2)$, our result indicates the existence for all $a,b>0$ and $\beta>0$.

\end{abstract}

\maketitle

\s{Introduction}
\renewcommand{\theequation}{1.\arabic{equation}}

Time-dependent systems of coupled nonlinear Schr\"odinger equations of the form
\beq\lab{eq:T-coupled-system}
\begin{cases}
-i\frac{\partial}{\partial t}\Phi_1=\Delta \Phi_1+g_1(|\Phi_1|^2)\Phi_1+\partial_1\varphi(|\Phi_1|^2,|\Phi_2|^2)\Phi_1,\\
-i\frac{\partial}{\partial t}\Phi_2=\Delta \Phi_2+g_2(|\Phi_2|^2)\Phi_2+\partial_2\varphi(|\Phi_1|^2,|\Phi_2|^2)\Phi_2,\\
\Phi_j=\Phi_j(x,t)\in \C, \, j=1,2,N\geq 1,
\end{cases}\quad (x,t)\in \R^N\times \R,
\eeq
are used as model for various physical phenomena, for instance binary mixtures of Bose-Einstein condensates, or the propagation of mutually incoherent wave packets in nonlinear optics; see  e.g.\ \cite{AkhmedievAnkiewicz.1999, Esry1998, Frantzeskakis2010, Timmermans1998}.
The ansatz $\Phi_j(x,t)=e^{i\lambda_jt}u_j(x), j=1,2$ for solitary wave solutions leads to the elliptic system
\beq\lab{eq:elliptic-system}
\begin{cases}
-\Delta u_1+\lambda_1u_1=f_1(u_1)+\partial_1 H(u_1, u_2),\\
-\Delta u_2+\lambda_2 u_2=f_2(u_2)+\partial_2 H(u_1,u_2),
\end{cases}\quad \hbox{in}\;\R^N,
\eeq
with $f_j(u_j)=g_j(|u_j|^2)u_j, j=1,2$ and $H(u_1,u_2)=\frac{1}{2}\varphi(|u_1|^2, |u_2|^2)$.

Since the masses
$$\int_{\R^N}|\Phi_j|^2\mathrm{d}x, j=1,2$$
are preserved along trajectories of \eqref{eq:T-coupled-system}, it is natural to consider them as prescribed. A natural approach to obtain solutions to \eqref{eq:elliptic-system} satisfying the normalization constraints
\begin{equation}\label{eq:norm}
 \int_{\R^N}|u_j|^2=a_j, j=1,2
\end{equation}
consists in finding critical points $(u_1,u_2)\in H^1(\R^N,\R^2)$ of the energy functional
$$\frac{1}{2}\int_{\R^N}\left(|\nabla u_1|^2+|\nabla u_2|^2\right) \mathrm{d}x-\int_{\R^N}\left(F_1(u_1)+F_2(u_2)+H(u_1,u_2)\right) \mathrm{d}x$$
 under the constraint \eqref{eq:norm}, where $F_j(s):=\int_0^s f_j(\tau)d\tau$ is the primitive function of $f_j, j=1,2$. Then the parameters $\lambda_1,\lambda_2$ appear as Lagrange multipliers.
\medskip


In the present paper we focus on the existence of solutions $(\lambda_1,\lambda_2, u, v)\in \R^2\times H^1(\R^N,\R^2)$ to the specific problem
\beq\lab{eq:20220902-maine1}
\begin{cases}
-\Delta u+\lambda_1u=\mu_1 |u|^{p-2}u+\beta r_1|u|^{r_1-2}u|v|^{r_2}\quad &\text{in $\R^N$,}\\
-\Delta v+\lambda_2v=\mu_2 |v|^{q-2}v+\beta r_2|u|^{r_1}|v|^{r_2-2}v\quad&\text{in $\R^N$,}\\
\int_{\R^N}|u|^2 \mathrm{d}x=a, \, \int_{\R^N}|v|^2 \mathrm{d}x=b.
\end{cases}
\eeq
Here $\mu_1,\mu_2,\beta>0$ are prescribed as are the masses $a,b>0$. The corresponding energy functional $J: \mathcal{H}\mapsto \R$ is
defined by
\begin{equation*}
J[u,v]=\frac{1}{2}[\|\nabla u\|_2^2+\|\nabla v\|_2^2]-\frac{\mu_1}{p}\|u\|_p^p-\frac{\mu_2}{q}\|v\|_q^q-\beta \int_{\R^N} |u|^{r_1} |v|^{r_2} \mathrm{d}x,
\end{equation*}
where
$\mathcal{H}:=H^1(\R^N)\times H^1(\R^N)$ is equipped with the norm
$$\|[u,v]\|_{\mathcal{H}}:=\left(\|u\|_{H^1}^{2}+\|v\|_{H^1}^{2}\right)^{\frac{1}{2}}.$$
The normalization constraint being $S_a \times S_b$ where
$$ S_a:=\left\{u\in H^1(\R^N): \|u\|_2^2=a\right\} \quad \mbox{and} \quad S_b:=\left\{u\in H^1(\R^N): \|u\|_2^2=b\right\}.$$
It is standard to show that $J: \mathcal{H}\mapsto \R$ is of class $C^1$ and that $J'$ takes bounded sets to bounded sets. \smallskip

For the case $N=3, r_1=r_2=2,\mu_1=\mu_2, p=q=4$, assuming that $\beta<0$, the existence of infinitely many positive solutions (namely $u>0$ and $v>0$) to
Problem \eqref{eq:20220902-maine1} was obtained by Bartsch and Soave in \cite{Bartsch2019}. When $\beta>0$,
 Gou and Jeanjean in \cite{Gou2018}  obtain a multiplicity result.  More precisely, under one of the following conditions
\begin{itemize}
\item[$(H_0)$] $N\geq 1, 1<p,q<2+\frac{4}{N}, r_1,r_2>1, 2+\frac{4}{N}<r_1+r_2<2^*$; \smallskip

\item[$(H_1)$] $N\geq 1, 2+\frac{4}{N}<p,q<2^*, r_1,r_2>1, r_1+r_2<2+\frac{4}{N}$,
\end{itemize}
two couples of solutions to Problem \eqref{eq:20220902-maine1} are obtained. One is a local minimizer for $J$ constrained on $S_a\times S_b$, the other is obtained through a constrained mountain pass or a constrained linking. In both cases it is assumed that $0<\beta<\beta_0$ is small enough.\smallskip

For the pure mass super-critical case, i.e., $2+\frac{4}{N}<p,q,r_1+r_2<2^*$, Bartsch, Jeanjean and Soave in \cite{Bartsch2016} consider the case of $N=3, p=q=4, r_1=r_2=2$. They obtain the existence of positive solutions to Problem \eqref{eq:20220902-maine1} provided $0<\beta<\beta_1(a,b)$ or $\beta>\beta_2(a,b)$. A more general case is considered by  Bartsch and Jeanjean in \cite{Bartsch2018} but it is also required there that $0<\beta<\beta_1(a,b)$  or $\beta>\beta_2(a,b)$. The values $\beta_1(a,b)$ and $\beta_2(a,b)$ depend on $a,b$ heavily. Typically, $\beta_1(a,b)\rightarrow 0$ as $a/b\rightarrow 0$ or $a/b \rightarrow +\infty$, while $\beta_{2}(a,b)\rightarrow +\infty$ as $a/b \rightarrow 0$ or $a/b \rightarrow +\infty$. In particular, there is no value of $\beta >0$ so that the results of \cite{Bartsch2016,Bartsch2018} yield a solution to Problem \eqref{eq:20220902-maine1} for all masses. \smallskip

In view of the above considerations, it is natural to investigate whether there exists some $\beta >0$ for which the existence of solutions is valid for all masses.  Actually, this issue was already raised as an {\it open problem} by Bartsch, Jeanjean and Soave, see \cite[Remark 1.3-(a)and (d)]{Bartsch2016}. In \cite{Bartsch2021}, Bartsch, Zhong and Zou  gave a first answer.
Precisely, assuming that $N=3, p=q=4, r_1=r_2=2$, they develop a new  approach based on fixed point index in cones, bifurcation theory, and a continuation method to deal with Problem \eqref{eq:20220902-maine1}. For the range of $\beta\in (0,\tau_0\min\{\mu_1,\mu_2\}]\cup(\tau_0\max\{\mu_1,\mu_2\},+\infty)$, where $\tau_0 \in (0,1)$ is explicitly known, the existence result holds for all masses. For the range of $\beta\in (\tau_0\min\{\mu_1,\mu_2\},\tau_0\max\{\mu_1,\mu_2\}]$  whether the existence result holds depends on the masses $a$ and $b$, more precisely, on the ratio $a/b$. For example, when $\mu_2<\mu_1$, for $a/b$ small, the existence result is satisfied. On the contrary when $\mu_2\leq \beta\leq \tau_0\mu_1$, there is no positive normalized solution provided $a/b$ is large. \smallskip

Our aim in the present paper is to make further advances on Bartsch-Jeanjean-Soave's open problem. Namely to explicit new ranges of $\beta >0$ for which there exist a solution to Problem \eqref{eq:20220902-maine1} for any $a>0$ and $b>0$. To appreciate our results observe that in \cite{Bartsch2021} all the nonlinear terms share the same homogeneity. In our results $p,q,r_1+r_2$ are not necessary the same.

\br\lab{remark:20230214-r1}
For the fixed frequencies problem (that is when $\lambda_1,\lambda_2$ are given) when $N=3, r_1=r_2=2,\mu_1=\mu_2, p=q=4$, usually one also need some structural conditions on the parameters to guarantee the existence of positive solutions. In that case there is also an open problem, which is called Sirakov's open problem in the literature, concerning the determination of the best range for the existence of positive solutions. In that direction we refer to \cite[Remark 4]{Sirakov2007} and \cite{WeiZhongZou-2022}.
\er

\s{Hypotheses and statement of results}\lab{sec:Hypotheses-statements}
\renewcommand{\theequation}{2.\arabic{equation}}

\noindent
In order to state our main results we first recall some classical facts. \smallskip

For $p\in (2,2^*)$, it is well known that the following scalar problem
\beq\lab{eq:def-Up}
  \begin{cases}
  -\Delta u+u=u^{p-1}\quad \;&\hbox{in}\;\R^N,\\
  u>0\quad &\hbox{in}~\R^N,\\
  u(0)=\max_{x\in \R^N} u(x) ~~\hbox{and}~~u\in H^1(\R^N)
  \end{cases}
\eeq
has a unique solution, denoted by $U_p$, which is a radial function, cf.\ \cite{Kwong1989}. \medskip

For $a,\mu\in \R^+$ fixed, consider the solutions $(\lambda, \omega)\in \R\times H^1(\R^N)$ with $\lambda>0$, of
\beq\lab{eq:w-unique}
\begin{cases}
-\Delta \omega+\lambda \omega=\mu \,  \omega^{p-1} \;&\hbox{in}\;\R^N,\\
w>0\quad &\hbox{in}~\R^N,\\
\omega(0)=\max_{x\in \R^N} \omega(x) ~~\hbox{and}~~u\in H^1(\R^N), \\
\|\omega\|_2^2=a.
\end{cases}
\eeq
It is also classical that these solutions are unique and can be obtained by scaling the solution to \eqref{eq:def-Up}.
Precisely, for $p\neq 2+\frac{4}{N}$, the couple
\beq\lab{eq:20220905-ze1}
\lambda=\lambda_{p,\mu,a}:=\mu^{-\frac{4}{N(p-2)-4}}\|U_p\|_{2}^{\frac{4(p-2)}{N(p-2)-4}} a^{-\frac{2(p-2)}{N(p-2)-4}}>0
\eeq
and $$U_{ p,\mu,\lambda}(x):=\left(\frac{\lambda}{\mu}\right)^{\frac{1}{p-2}} U_p(\sqrt{\lambda}x),$$  solves \eqref{eq:w-unique}.
We denote this unique solution by
\beq\lab{def:zpua}
z_{p,\mu,a}(x):= U_{p,\mu,\lambda_{p,\mu,a}}(x).
\eeq


\begin{definition}\lab{def:20220930-d1}
Let $a >0$, $b >0$ be given. We call $[u,v]$ a ground state solution (or energy ground state) to Problem \eqref{eq:20220902-maine1} if $[u,v]$ is a solution to Problem \eqref{eq:20220902-maine1} and
\begin{equation*}
J[u,v]=\inf\left\{J[\phi,\psi]: [\phi,\psi]~\hbox{is a solution to Problem \eqref{eq:20220902-maine1}}\right\}.
\end{equation*}
\end{definition}
\smallskip

Finally, let $p,q,\mu_1,\mu_2$ be fixed and for $a>0$ define
\begin{align*}
b_{p,q,\mu_1,\mu_2,a}:=&\left(\frac{\mu_2}{\mu_1} \frac{(q-2)N-4}{(p-2)N-4} \frac{2N-(N-2)p}{2N-(N-2)q}\right)^{\frac{(q-2)N-4}{2N-(N-2)q}}
\|U_q\|_{2}^{\frac{4(q-2)}{2N-(N-2)q}}  \nonumber\\
& \times \|U_p\|_{2}^{-\frac{4(p-2)}{2N-(N-2)q} \frac{(q-2)N-4}{(p-2)N-4}} a^{\frac{2N-(N-2)p}{2N-(N-2)q} \frac{(q-2)N-4}{(p-2)N-4}}.
\end{align*}
and
\beq\lab{def:-best-const}
\beta_{p,\mu,a,N,r}:=\frac{1}{2}\inf_{h\in H^1(\R^N)\backslash\{0\}}\frac{\|\nabla h\|_2^2}{\int_{\R^N}|z_{p,\mu,a}|^{r} |h|^2 \mathrm{d}x}.
\eeq
\smallskip
We remark that $J[0,z_{q,\mu_2,b}]<J[z_{p,\mu_1,a},0]$ if and only if $b>b_{p,q,\mu_1,\mu_2,a}$ (by a direct calculation from Lemma \ref{lemma:20210622-l2}-(i)). Now we can state our main result.
\bt\lab{th:main-t1}
Let $1\leq N\leq 4, 2+\frac{4}{N}<p,q,r_1+r_2<2^*$ and $r_1,r_2>1$. Then the following holds true.
\begin{itemize}
\item[(i)] For $a>0$ and $b\in [b_{p,q,\mu_1,\mu_2,a},+\infty)$. If
$$\begin{cases}
r_1<2\\
\beta>0
\end{cases}~\hbox{or}~
\begin{cases}
r_1=2\\
\beta>\beta_{q,\mu_2,b,N,r_2} \, ,
\end{cases}$$
then there exists a ground state solution $(\lambda_1,\lambda_2,u,v)$ to Problem \eqref{eq:20220902-maine1}. In addition, $\lambda_1>0,\lambda_2>0$ and $u, v $ are Schwartz symmetric functions.
\item[(ii)]Let $a>0$ and $b\in (0,b_{p,q,\mu_1,\mu_2,a}]$. If
$$\begin{cases}
r_2<2\\
\beta>0
\end{cases}~ \hbox{or}~
\begin{cases}
r_2=2\\
\beta>\beta_{p,\mu_1,a,N,r_1} \, ,
\end{cases}$$
then there exists a ground state solution $(\lambda_1,\lambda_2,u,v)$ to Problem \eqref{eq:20220902-maine1}. In addition, $\lambda_1>0,\lambda_2>0$ and $u, v $ are Schwartz symmetric functions.
\end{itemize}
\et

In Lemmas \ref{lemma:20220907-l1} and \ref{20220905-l1}, the behavior of $\beta_{p,\mu,a,N,r}$ with respect to the various parameters is studied. In particular, for the case of $N=1,2$, we prove that $\beta_{p,\mu,a,N,r}=0$ always is true, so we  have the following conclusion.
\bc\lab{cro:main-c2}
For  $N=1,2$, assume that $2+\frac{4}{N}<p,q$ and $r_1>1,r_2>1, r_1+r_2>2+\frac{4}{N}$. Then
\begin{itemize}
\item[(i)] for any $a>0$ and $\beta>0$, there exists a ground state solution $(\lambda_1,\lambda_2,u,v)$ to Problem \eqref{eq:20220902-maine1} provided $b\in[b_{p,q,\mu_1,\mu_2,a},+\infty)$ and $r_1\leq 2$.  In addition, $\lambda_1>0,\lambda_2>0$ and $u, v $ are Schwartz symmetric functions.
\item[(ii)]for any $a>0$ and $\beta>0$, there exists a positive ground state solution $(\lambda_1,\lambda_2,u,v)$ to Problem \eqref{eq:20220902-maine1} provided  $b\in(0,b_{p,q,\mu_1,\mu_2,a}]$ and $r_2\leq 2$. In addition, $\lambda_1>0,\lambda_2>0$ and $u, v $ are Schwartz symmetric functions.
\end{itemize}
\ec

For $N=3,4$, we get directly from Theorem \ref{th:main-t1}.
\bc\lab{cro:main-c1}
For $N=3$ or $N=4$, assume that $2+\frac{4}{N}<p,q, r_1+r_2<2^*$ and $1<r_1,r_2<2$. Then for any $a,b,\beta>0$, there exists a ground state solution $(\lambda_1,\lambda_2,u,v)$ to Problem \eqref{eq:20220902-maine1}. In addition, $\lambda_1>0,\lambda_2>0$ and $u, v $ are Schwartz symmetric functions.
\ec

Let us now describe the main steps of the proof of Theorem \ref{th:main-t1}. \smallskip

The ground states will be obtain as minimum of the functional $J$ on the natural constraint
\begin{equation}
\lab{eq:lla}
\mathcal{P}:=\left\{[u,v]\in \mathcal{H}\backslash\{[0,0]\}: P[u,v]=0\right\},
\end{equation}
where
\begin{align*}
P[u,v]:=&\|\nabla u\|_2^2+\|\nabla v\|_2^2-\frac{(p-2)N}{2p}\mu_1\|u\|_p^p-\frac{(q-2)N}{2q}\mu_2\|v\|_q^q\nonumber\\
&-\frac{(r_1+r_2-2)N}{2}\beta \int_{\R^N} |u|^{r_1} |v|^{r_2} \mathrm{d}x.
\end{align*}
Here also one can show, in a standard way, that $P$ is of class $C^1$ and that $P'$ takes bounded sets into bounded sets.
This constraint is a natural one in the sense that, as a consequence of a Pohozaev identity, all critical points of $J$ constrained to $S_a \times S_b$ belong to $\mathcal{P}$. Thus, if we manage to show that
$$\inf_{\mathcal{P} \cap (S_a \times S_b)}J[u,v]$$
is achieved by a critical point this will prove the existence of a ground state. \medskip

Instead of attacking directly this problem we consider a relaxed one.
For $a\geq 0$, we introduce the set
\begin{equation*}
D_a:=\{u\in H^1(\R^N): \|u\|_2^2\leq a\}
\end{equation*}
and for any $(a,b)\in \R^+\times \R^+$, we define
\begin{equation*}
\mathcal{P}_{(a,b)}:=\mathcal{P}\cap (D_a\times D_b).
\end{equation*}

The relaxed problem is :
\begin{align}
\lab{eq:20220930-e1}
\mbox{to find a critical point } & [u,v] \in \mathcal{H} \backslash \{[0,0]\} \mbox{ of } J  \mbox{ constrained to } D_a \times D_b\nonumber\\
& \mbox{ at the level }
C_{(a,b)}:=\inf_{\mathcal{P}_{(a,b)}}J[u,v].
\end{align}


Once this problem is settle, we show, see Lemma \ref{lemma:20220831-l1}, that $[u,v]\in S_a\times S_b$ solving thus the initial problem. \medskip

This idea of first considering a relaxed problem is reminiscent of the works \cite{BM21,MeSc2022}.  However the strategy of the proofs in \cite{BM21,MeSc2022} and in the present work are essentially distinct.

To solve Problem \eqref{eq:20220930-e1}, we shall make use of a classical result due to Ghoussoub. This result \cite[Theorem 4.1]{Gh} provides a  sequence $\{(u_n, v_n)\} \subset \mathcal{H}$ for $J$ restricted to $D_a \times D_b$ at the level $C_{(a,b)}$ which have useful additional properties. In particular it is {\it localized} around $\mathcal{P}$ and this insure its boundedness. \smallskip

 To obtain the convergence of $\{(u_n, v_n)\} \subset \mathcal{H}$, and thus the existence of a critical point at this level, a key point is to show that the level $C_{(a,b)}$ is {\it sufficiently low}. More precisely, observe that
 $[z_{p,\mu_1,a},0], [0,z_{q,\mu_2,b}]\in \mathcal{P}_{(a,b)}$, and thus that it holds
\begin{equation*}
C_{(a,b)}\leq \min\{m_{p,\mu_1,a}, m_{q,\mu_2,b}\},
\end{equation*}
where we have set
$$ m_{p,\mu_1,a}:= J(z_{p,\mu_1,a}, 0) \quad \mbox{and} \quad m_{q,\mu_2,b}:= J(0, z_{q,\mu_2,b}).$$
The condition that we use to conclude to the convergence of $\{(u_n, v_n)\} \subset \mathcal{H}$ is
\beq\lab{eq:20221001-xe2}
C_{(a,b)}< \min\{m_{p,\mu_1,a}, m_{q,\mu_2,b}\}.
\eeq
To prove \eqref{eq:20221001-xe2} we study the behavior of $J$ restricted to $\mathcal{P}_{(a,b)}$ locally around $[z_{p,\mu_1,a},0]$ and $[0,z_{q,\mu_2,b}]$ in a way which seems new to us. In particular the inequality \eqref{eq:20221001-xe2} is proved testing $J$ on a point in $ \mathcal{P}_{(a,b)} \backslash \{S_a \times S_b\}$.
It is at this step, and only there, that appears the need to restrict our results to small values of $r_1 >0$ and $r_2 >0$. Having proved \eqref{eq:20221001-xe2} the convergence of $\{(u_n, v_n)\} \subset \mathcal{H}$ follows, leading to the existence of a ground state. \smallskip

To obtain the additional properties of the ground states stated in Theorem \ref{th:main-t1}, and its corollaries, we use the fact that any minimum for
$$\inf_{\mathcal{P} \cap (S_a \times S_b)}J[u,v]$$
is a critical point for $J$ constrained to $S_a \times S_b$ and that it is achieved by Schwartz symmetric functions. \smallskip

The paper is organized as follows. In Section \ref{sec:prelim} we study the scalar problems which appear in the discussion of the convergence of our Palais-Smale sequence $\{(u_n, v_n)\} \subset \mathcal{H}$.
 Section \ref{sec:Pohozaev} is devoted to the analysis of the geometry of the set $\mathcal{P}_{(a,b)}$. We also prove that any critical point for $J$ restricted to $\mathcal{P}_{(a,b)}$ is a critical point for $J$ restricted to $D_a \times D_b$, namely that no Lagrange multiplier arises in the presence of the constraint $\mathcal{P}$. In
Section \ref{sec:Rearrangement} we show that it is not restrictive to search ground states within the Schwartz symmetric functions. In Section \ref{sec:Palais_Smale} we establish the existence of our particular Palais-Smale sequence $\{[u_n,v_n]\} \subset \mathcal{H}$. Section \ref{sec:estimation} is fully devoted to the proof of the strict inequality \eqref{eq:20221001-xe2}. Finally, in Section \ref{sec:PS-sequence}, we establish the convergence of  $\{[u_n,v_n]\} \subset \mathcal{H}$  and we give the proof of Theorem \ref{th:main-t1}.



\s{Some  Preliminaries}\lab{sec:prelim}
\renewcommand{\theequation}{3.\arabic{equation}}
We present in this section results which are slightly more general than needed. \smallskip

Let us introduce the set
\begin{equation*}
\tilde{\mathcal{P}}_{p,\mu}:=\left\{u\neq 0: \|\nabla u\|_2^2=\frac{(p-2)N}{2p}\mu\|u\|_p^p\right\}
\end{equation*}
and for $a>0$,
\begin{equation*}
\tilde{\mathcal{P}}_{p,\mu,a}:=\tilde{\mathcal{P}}_{p,\mu}\cap S_a=\left\{\|u\|_2^2=a: \|\nabla u\|_2^2=\frac{(p-2)N}{2p}\mu\|u\|_p^p\right\}.
\end{equation*}

Recalling the fiber map
\beq\lab{eq:ut-invariant}
u(x)\mapsto (t\star u)(x):=t^{\frac{N}{2}}u(tx),
\eeq
for $(t, u)\in \R^+\times H^1(\R^N)$, which preserves the $L^2$-norm. We have the following results.
\bl\lab{lemma:20210622-l1}
\begin{itemize}
\item[(i)]If $u$ is a solution to \eqref{eq:w-unique}, then $u\in \tilde{\mathcal{P}}_{p,\mu,a}$.  Also the unique solution to \eqref{eq:w-unique} minimizes $I_{p,\mu}$ on $\tilde{\mathcal{P}}_{p,\mu,a}$.
\item[(ii)]If $p\neq 2+\frac{4}{N}$, then for any $u\neq 0$, there exists a unique $t=t(u)>0$ such that $t\star u\in \tilde{\mathcal{P}}_{p,\mu}$ and it holds that
\begin{equation}
\lab{eq:20210622-xe1}
I_{p,\mu}[t\star u]=\begin{cases}
\max_{s>0} I_{p,\mu}[s\star u]\quad &\hbox{if}~2+\frac{4}{N}<p<2^*,\\
\min_{s>0} I_{p,\mu}[s\star u]\quad &\hbox{if}~2<p<2+\frac{4}{N}.
\end{cases}
\end{equation}
\item[(iii)]Let $z_{p,\mu,a}$ be defined by \eqref{def:zpua}. Then
\begin{align*}
m_{p,\mu,a}:=&I_{p,\mu}[z_{p,\mu,a}]=\inf_{u\in \tilde{\mathcal{P}}_{p,\mu,a}}I_{p,\mu}[u]\\
=&\begin{cases}
    \inf_{u\in S_a}\max_{s>0} I_{p,\mu}[s\star u]\quad &\hbox{if}~2+\frac{4}{N}<p<2^*,\\
    \inf_{u\in S_a}\min_{s>0} I_{p,\mu}[s\star u]=\inf_{u\in S_a}I_{p,\mu}[u]\quad &\hbox{if}~2<p<2+\frac{4}{N}.
    \end{cases}
\end{align*}
\end{itemize}
\el
\bp
(i) We refer to \cite[Lemma 2.1]{Bartsch2016}.\\
(ii) For $u\neq 0$ and $s>0$, a direct calculation shows that
\begin{equation*}
f(s):=I_{p,\mu}[s\star u]=\frac{1}{2}\|\nabla u\|_2^2 s^2-\frac{\mu}{p}\|u\|_p^p s^{\frac{(p-2)N}{2}}.
\end{equation*}
It is easy to see that $s\star u \in \tilde{\mathcal{P}}_{p,\mu}$ if and only if $s$ is a critical point of $f(s)$.
Then if $p\neq 2+\frac{4}{N}$, one can show that $f(s)$ has an unique critical point $$t=t(u):=\left(\frac{\|\nabla u\|_2^2}{\frac{(p-2)N}{2p}\mu \|u\|_p^p}\right)^{\frac{2}{(p-2)N-4}}$$
and
$$f(t)=\frac{(p-2)N-4}{4p} \left(\frac{2p}{(p-2)N}\right)^{\frac{(p-2)N}{(p-2)N-4}} \left(\|\nabla u\|_2^2\right)^{\frac{(p-2)N}{(p-2)N-4}} \left(\mu\|u\|_p^p\right)^{\frac{4}{(p-2)N-4}}.$$
Furthermore, $f(t)>0$ achieves the maximum if $p>2+\frac{4}{N}$, while $f(t)<0$ achieves the minimum if $p<2+\frac{4}{N}$. Hence, \eqref{eq:20210622-xe1} holds.\\
(iii) By (i) and (ii), it is easy to see.
\ep

\bl\lab{lemma:20210622-l2}
\begin{itemize}
\item[(i)]
If $2+\frac{4}{N}<p<2^*$, then $m_{p,\mu,a}$ is continuous and  decreases strictly respect to $a\in \R^+$, furthermore,
$$\lim_{a\rightarrow 0^+}m_{p,\mu,a}=+\infty \;\hbox{and}\;\lim_{a\rightarrow +\infty}m_{p,\mu,a}=0.$$
\item[(ii)]If $2<p<2+\frac{4}{N}$, then $m_{p,\mu,a}$ is continuous and  decreases strictly with respect to $a\in \R^+$, furthermore,
$$\lim_{a\rightarrow 0^+}m_{p,\mu,a}=0 \;\hbox{and}\;\lim_{a\rightarrow +\infty}m_{p,\mu,a}=-\infty.$$
\end{itemize}
\el
\bp
By Lemma \ref{lemma:20210622-l1}-(iii), $m_{p,\mu,a}=I_{p,\mu}[z_{p,\mu,a}]$. A direct calculation shows that
 \begin{equation*}
 m_{p,\mu,a}=\frac{1}{2}\frac{(p-2)N-4}{2N-(N-2)p} \left(\|U_p\|_2^2\right)^{\frac{2(p-2)}{(p-2)N-4}} \mu a^{-\frac{2N-(N-2)p}{(p-2)N-4}}.
 \end{equation*}
  If $2+\frac{4}{N}<p<2^*$, we see the exponent $-\frac{2N-(N-2)p}{(p-2)N-4}<0$, hence $m_{p,\mu,a}$ decreases strictly  with respect to $a\in \R^+$. Furthermore,
 $$\lim_{a\rightarrow 0^+}m_{p,\mu,a}=+\infty \;\hbox{and}\;\lim_{a\rightarrow +\infty}m_{p,\mu,a}=0.$$
 Similarly, if $2<p<2+\frac{4}{N}$, then $-\frac{2N-(N-2)p}{(p-2)N-4}>0$, and noting that
 $\frac{1}{2}\frac{(p-2)N-4}{2N-(N-2)p}<0$ in such a case,
 thus $m_{p,\mu,a}$ also decreases strictly respect to $a\in \R^+$. Furthermore,
 $$\lim_{a\rightarrow 0^+}m_{p,\mu,a}=0 \;\hbox{and}\;\lim_{a\rightarrow +\infty}m_{p,\mu,a}=-\infty.$$
\ep


\s{A Pohozaev type constraint and its properties}\lab{sec:Pohozaev}
\renewcommand{\theequation}{4.\arabic{equation}}

In this section we establish various properties of the constraint $\mathcal{P}$ and we prove, in particular, it does not generate a Lagrange parameter. First observe that,
\bl\lab{lemma:20221001-l1}
Let $\mathcal{P}$ be the so-called Pohozaev manifold defined by \eqref{eq:lla}.
If  $[u,v]\in \mathcal{H}\backslash\{[0,0]\}$ is a solution to
\begin{equation}
\lab{eq:main-equation}
\begin{cases}
-\Delta u+\lambda_1u=\mu_1 |u|^{p-2}u+\beta r_1|u|^{r_1-2}u|v|^{r_2}\quad &\text{in $\R^N$,}\\
-\Delta v+\lambda_2v=\mu_2 |v|^{q-2}v+\beta r_2|u|^{r_1}|v|^{r_2-2}v\quad&\text{in $\R^N$,}
\end{cases}
\end{equation}
 then $[u,v]\in \mathcal{P}$.
\el
\bp
This result can be shown as in \cite[Lemma 4.6]{Bartsch2016}.
\ep

Now, for any fixed $[u,v]\in\mathcal{H}\backslash\{[0,0]\}$, let us define $\Psi_{[u,v]}:\R^+\rightarrow \R$ by
\begin{align}
\Psi_{[u,v]}(t):=J[t\star u, t\star v]=&\frac{1}{2}\left[\|\nabla u\|_2^2+\|\nabla v\|_2^2\right]t^2
-\frac{\mu_1}{p}\|u\|_p^p  \, t^{\frac{(p-2)N}{2}}\nonumber\\
&-\frac{\mu_2}{q}\|v\|_q^q \, t^{\frac{(q-2)N}{2}} -\beta \Big(\int_{\R^N} |u|^{r_1} |v|^{r_2} \mathrm{d}x \Big) \,  t^{\frac{(r_1+r_2-2)N}{2}} \nonumber.
\end{align}
A direct calculation shows that
\begin{align}
\lab{eq:bu-20210623-e1}
\Psi'_{[u,v]}(t)=&\frac{d}{dt} J[t\star u, t\star v]\nonumber\\
=&\left[\|\nabla u\|_2^2+\|\nabla v\|_2^2\right]t-\frac{(p-2)N}{2p}\mu_1\|u\|_p^p \, t^{\frac{(p-2)N-2}{2}}\nonumber\\
&-\frac{(q-2)N}{2q}\mu_2\|v\|_q^q \, t^{\frac{(q-2)N-2}{2}}-\frac{(r_1+r_2-2)N}{2}\beta \Big(\int_{\R^N} |u|^{r_1} |v|^{r_2} \mathrm{d}x \Big) \,  t^{\frac{(r_1+r_2-2)N-2}{2}} .
\end{align}
We remark that $P[u,v]=\Psi'_{[u,v]}(1)$ and $P[t\star u, t\star v]=t \Psi'_{[u,v]}(t)$.

\bl\lab{lemma:20210622-wbul2}
 Let $[u,v]\in\mathcal{H}\backslash\{[0,0]\}$. Then: $t\in \R^+$ is a critical point of $\Psi_{[u,v]}(t)$ if and only if $[t\star u, t\star v]\in \mathcal{P}$.
\el
\bp
By the fact $P[t\star u, t\star v]=t \Psi'_{[u,v]}(t)$ and $t\neq 0$, we see that
$$[t\star u, t\star v]\in \mathcal{P}\Leftrightarrow P[t\star u, t\star v]=0\Leftrightarrow \Psi'_{[u,v]}(t)=0.$$
\ep

\bc\lab{cro:uniqueness-t}
Assume that $2+\frac{4}{N}<p,q,r_1+r_2<2^*$.
Then for any $[0,0]\neq [u,v]\in D_a\times D_b$, there exists a unique $t=t_{[u,v]}>0$ such that $[t\star u, t\star v]\in \mathcal{P}_{(a,b)}$. Furthermore, $t_{[u,v]}<(resp.~=,>) 1$ if and only if $P[u,v]<(resp.~=,>) 0$.
\ec
\bp
Noting that the fiber map \eqref{eq:ut-invariant} preserves the $L^2$-norm, we have that $[t\star u, t\star v]\in D_a\times D_b$ for all $t>0$. Then by Lemma \ref{lemma:20210622-wbul2}, $[t\star u, t\star v]\in \mathcal{P}_{(a,b)}$ if and only if $\Psi'_{[u,v]}(t)=0$. Now by \eqref{eq:bu-20210623-e1} and $2+\frac{4}{N}<p,q,r_1+r_2<2^*$, we obtain the uniqueness of $t=t_{[u,v]}$. Furthermore, $t_{[u,v]}$ attains the maximum of $\Psi_{[u,v]}(t)$ in $t\in\R^+$ and
$$\Psi'_{[u,v]}(t)>0\;\hbox{for}~0<t<t_{[u,v]}\;\hbox{while}~\Psi'_{[u,v]}(t)<0\;\hbox{for}~t>t_{[u,v]}.$$
So combining with $P[u,v]=\Psi'_{[u,v]}(1)$, we obtain that
$$P[u,v]<(resp.~=,>)0\Leftrightarrow \Psi'_{[u,v]}(1)<(resp.~=,>)0\Leftrightarrow t_{[u,v]}<(resp.~=,>)1.$$
\ep

\bl\lab{lemma:20210623-xbzzl1}
Suppose that $2+\frac{4}{N}<p,q,r_1+r_2<2^*$.
There exists some $C_0>0$  depending only on $p,q,r_1,r_2,N$,  such that
\begin{equation*}
J[u,v]\geq C_0 \Big[\|\nabla u\|_2^2+\|\nabla v\|_2^2\Big], \quad \mbox{for all }  [u,v]~\hbox{satisfying}~P[u,v]=0.
\end{equation*}
\el
\bp
For $2+\frac{4}{N}<p,q,r_1+r_2<2^*$,
set
\begin{equation*}
\tau:=\max\left\{\frac{2}{(p-2)N}, \frac{2}{(q-2)N}, \frac{2}{(r_1+r_2-2)N}\right\},
\end{equation*}
then one can see that $0<\tau<\frac{1}{2}$.
By the definition of $\tau$,
\begin{align}
&\frac{\mu_1}{p}\|u\|_p^p+\frac{\mu_2}{q}\|v\|_q^q+\beta \int_{\R^N} |u|^{r_1}|v|^{r_2}\mathrm{d}x\nonumber\\
\leq &\tau \left\{\frac{(p-2)N}{2p}\mu_1\|u\|_p^p+\frac{(q-2)N}{2q}\mu_2\|v\|_q^q
+\frac{(r_1+r_2-2)N}{2}\beta \int_{\R^N} |u|^{r_1} |v|^{r_2} \mathrm{d}x\right\} \nonumber
\end{align}
So for any $[u,v]$ with $P[u,v]=0$,  we have that
\begin{equation*}
\frac{\mu_1}{p}\mu_1\|u\|_p^p+\frac{\mu_2}{q}\mu_2\|v\|_q^q+\beta \int_{\R^N} |u|^{r_1}|v|^{r_2}\mathrm{d}x
\leq \tau [\|\nabla u\|_2^2+\|\nabla v\|_2^2].
\end{equation*}
 Hence, we can take $C_0:=\frac{1}{2}-\tau>0$ such that
\begin{align*}
J[u,v]=&\frac{1}{2}[\|\nabla u\|_2^2+\|\nabla v\|_2^2]-\left[\frac{\mu_1}{p}\|u\|_p^p+\frac{\mu_2}{q}\|v\|_q^q+\beta \int_{\R^N} |u|^{r_1}|v|^{r_2}\mathrm{d}x\right]\\
\geq &\frac{1}{2}[\|\nabla u\|_2^2+\|\nabla v\|_2^2]-\tau [\|\nabla u\|_2^2+\|\nabla v\|_2^2]\\
=&C_0[\|\nabla u\|_2^2+\|\nabla v\|_2^2].
\end{align*}
\ep

\bc\lab{cro:20210623-c1}
Suppose that $2+\frac{4}{N}<p,q,r_1+r_2<2^*$. For any $a>0,b>0$,
$J\big|_{\mathcal{P}_{(a,b)}}$ is coercive, i.e.,
$$\lim_{\stackrel{[u,v]\in \mathcal{P}_{(a,b)}}{\|\nabla u\|_2^2+\|\nabla v\|_2^2\rightarrow +\infty}}J[u,v]=+\infty.$$
\ec
\bp
It follows by Lemma \ref{lemma:20210623-xbzzl1}.
\ep

\bl\lab{lemma:tidu-xyj}
Suppose that $2+\frac{4}{N}<p,q,r_1+r_2<2^*$.
For any given $a\geq 0,b\geq 0$ with $(a,b)\neq (0,0)$, there exists some $\delta_{(a,b)}>0$ such that
\begin{equation*}
\inf_{(u,v)\in \mathcal{P}_{(a,b)}} \Big[\|\nabla u\|_2^2+\|\nabla v\|_2^2 \Big]\geq \delta_{(a,b)}.
\end{equation*}
\el
\bp
Recalling the well known Gagliardo-Nirenberg inequality,
\begin{equation*}
\|u\|_p\leq C_{N,p} \|\nabla u\|_{2}^{\frac{N(p-2)}{2p}} \|u\|_{2}^{1-\frac{N(p-2)}{2p}}, \quad \mbox{for all } u\in H^1(\R^N) ~\hbox{and all }~ ~2<p<2^*,
\end{equation*}
for any $[u,v]\in \mathcal{P}$, $P[u,v]=0$ implies that at least one of the following holds.
\begin{itemize}
\item[(i)] $\displaystyle \frac{1}{3}[\|\nabla u\|_2^2+\|\nabla v\|_2^2]\leq \frac{(p-2)N}{2p}\mu_1\|u\|_p^p$,
\item[(ii)] $\displaystyle \frac{1}{3}[\|\nabla u\|_2^2+\|\nabla v\|_2^2]\leq \frac{(q-2)N}{2q}\mu_2\|v\|_q^q$,
\item[(iii)] $\displaystyle \frac{1}{3}[\|\nabla u\|_2^2+\|\nabla v\|_2^2]\leq \frac{(r_1+r_2-2)N}{2}\beta \int_{\R^N} |u|^{r_1} |v|^{r_2} \mathrm{d}x$.
\end{itemize}
If $(i)$ holds, we have that
\begin{align*}
\frac{1}{3}[\|\nabla u\|_2^2+\|\nabla v\|_2^2]\leq&  \frac{(p-2)N}{2p}\mu_1 C_{N,p}^{p}
\|\nabla u\|_{2}^{\frac{N(p-2)}{2}} \|u\|_{2}^{\frac{2N-(N-2)p}{4}}\\
\leq& \frac{(p-2)N}{2p}\mu_1 C_{N,p}^{p} \left(\|\nabla u\|_2^2+\|\nabla v\|_2^2\right)^{\frac{N(p-2)}{4}} \left(\|u\|_2^2+\|v\|_2^2\right)^{\frac{2N-(N-2)p}{4}}.
\end{align*}
By $2+\frac{4}{N}<p<2^*$, we have that $\frac{N(p-2)}{4}>1$ and $\frac{2N-(N-2)p}{4}>0$. Hence,
\begin{equation}
\lab{eq:20210623-pe1}
\|\nabla u\|_2^2+\|\nabla v\|_2^2\geq \left(\frac{3(p-2)N}{2p}\mu_1 C_{N,p}^{p}\right)^{-\frac{4}{N(p-2)-4}} \left(\|u\|_2^2+\|v\|_2^2\right)^{-\frac{2N-(N-2)p}{N(p-2)-4}}.
\end{equation}
Similarly, if $(ii)$ holds, we have
\begin{equation}
\lab{eq:20210623-pe2}
\|\nabla u\|_2^2+\|\nabla v\|_2^2\geq \left(\frac{3(q-2)N}{2q}\mu_2 C_{N,q}^{q}\right)^{-\frac{4}{N(q-2)-4}} \left(\|u\|_2^2+\|v\|_2^2\right)^{-\frac{2N-(N-2)q}{N(q-2)-4}}.
\end{equation}
If $(iii)$ holds, by H\"older inequality, a similar argument shows that
\beq\lab{eq:20210623-pe3}
\|\nabla u\|_2^2+\|\nabla v\|_2^2\geq \left(\frac{3(r_1+r_2-2)N}{2}\beta C_{N,r_1+r_2}^{r_1+r_2}\right)^{-\frac{4}{N(r_1+r_2-2)-4}} \left(\|u\|_2^2+\|v\|_2^2\right)^{-\frac{2N-(N-2)(r_1+r_2)}{N(r_1+r_2-2)-4}}.
\eeq
For $2+\frac{4}{N}<p,q,r_1+r_2<2^*$, by \eqref{eq:20210623-pe1}-\eqref{eq:20210623-pe3}, we can find some $\delta_{(a,b)}>0$ such that
\begin{equation*}
\|\nabla u\|_2^2+\|\nabla v\|_2^2\geq \delta_{(a,b)}, \quad \mbox{for all } [u,v]\in \mathcal{P}_{(a,b)}.
\end{equation*}
\ep

\bc\lab{cro:def-Cab}
Assume that $2+\frac{4}{N}<p,q,r_1+r_2<2^*$. For any $(a,b)\in \R^+\times \R^+$,
$$C_{(a,b)}=\inf_{[0,0]\neq [u,v]\in D_a\times D_b}\max_{t>0} J[t\star u, t\star v]>0.$$
\ec
\bp
For any $[u,v]\in \mathcal{P}_{(a,b)}$, by Lemma \ref{lemma:20210623-xbzzl1} and Lemma \ref{lemma:tidu-xyj}, we have that
$$J[u,v]\geq C_0[\|\nabla u\|_2^2+\|\nabla v\|_2^2]\geq C_0\delta_{(a,b)}>0.$$
Hence, $C_{(a,b)}$ is well defined and $C_{(a,b)}\geq C_0\delta_{(a,b)}>0$. Furthermore, by Corollary \ref{cro:uniqueness-t}, one can see that
 $$\inf_{\mathcal{P}_{(a,b)}}J[u,v]=\inf_{[0,0]\neq [u,v]\in D_a\times D_b}\max_{t>0} J[t\star u, t\star v].$$
\ep

We end this section by showing that any critical point of $J$ restricted to $\mathcal{P}_{(a,b)}$ is a critical point of $J$ restricted to $D_a \times D_b$. This implies, in particular, that any minimum of $C_{(a,b)}$  is a ground state to Problem \eqref{eq:20220902-maine1}.

\bl\lab{lemma:20210623-l1}
For any critical point of $J\big|_{\mathcal{P}_{(a,b)}}$, if $\Psi''_{[u,v]}(1)\neq 0$, then there exists some $\lambda_1,\lambda_2\in \R$ such that
\begin{equation*}
J'[u,v]+\lambda_1[u,0]+\lambda_2[0,v]=0.
\end{equation*}
\el

\bp
Firstly, it is easy to see that there exists some $\lambda_1,\lambda_2,\mu\in \R$ such that
\beq\lab{eq:20210623-xe1}
J'[u,v]+\lambda_1[u,0]+\lambda_2[0,v]+\mu P'[u,v]=0~\hbox{in $\mathcal{H}^{-1}$, the dual space of}~\mathcal{H}.
\eeq
Hence, we only need to prove that $\mu=0$. The functional associated to this
equation
is given by
\begin{equation*}
\Phi[u,v]:=J[u,v]+\frac{1}{2}\lambda_1\|u\|_2^2+\frac{1}{2}\lambda_2\|v\|_2^2+\mu P[u,v].
\end{equation*}
The fact that $[u,v]$ solves \eqref{eq:20210623-xe1} implies that $t=1$ is a critical point of
$\tilde{\Psi}_{[u,v]}:\R^+\rightarrow \R$ given by
\begin{align*}
\tilde{\Psi}_{[u,v]}(t) := \Phi[t\star u, t\star v]=&J[t\star u, t\star v]+\frac{\lambda_1}{2}\|u\|_2^2+\frac{\lambda_2}{2}\|v\|_2^2+\mu P[t\star u, t\star v]\nonumber\\
=&\Psi_{[u,v]}(t)+\frac{\lambda_1}{2}\|u\|_2^2+\frac{\lambda_2}{2}\|v\|_2^2+\mu t\Psi'_{[u,v]}(t).
\end{align*}
We have
\begin{equation*}
\tilde{\Psi}'_{[u,v]}(t) = \frac{d}{dt}\Phi[t\star u, t\star v]=(1+\mu)\Psi'_{[u,v]}(t)+\mu \Psi''_{[u,v]}(t).
\end{equation*}
Substituting $t=1$ and $\Psi'_{[u,v]}(1)=P[u,v]=0$, we obtain that
$\mu \Psi''_{[u,v]}(1)=0$. Then it follows $\mu=0$ by $\Psi''_{[u,v]}(1)\neq 0$.
\ep

\bl\lab{lemma:20210623-l2}
Suppose that $2+\frac{4}{N}<p,q,r_1+r_2<2^*$. Then for any $[u,v]\in \mathcal{P}$, we have $\Psi''_{[u,v]}(1)<0$.
\el
\bp
By a direct computation, we have that
\begin{align*}
\Psi''_{[u,v]}(t)=&[\|\nabla u\|_2^2+\|\nabla v\|_2^2]-\frac{(p-2)N}{2p}\frac{(p-2)N-2}{2}\mu_1\|u\|_p^p \, t^{\frac{(p-2)N-4}{2}}\nonumber\\
&-\frac{(q-2)N}{2q}\frac{(q-2)N-2}{2}\mu_2\|v\|_q^q \,  t^{\frac{(q-2)N-4}{2}}\nonumber\\
&-\frac{(r_1+r_2-2)N}{2}\frac{(r_1+r_2-2)N-2}{2}\beta \Big( \int_{\R^N} |u|^{r_1} |v|^{r_2} \mathrm{d}x \Big) t^{\frac{(r_1+r_2-2)N-4}{2}}.
\end{align*}
Hence,
\begin{align*}
\Psi''_{[u,v]}(1)=&[\|\nabla u\|_2^2+\|\nabla v\|_2^2]-\frac{(p-2)N}{2p}\frac{(p-2)N-2}{2}\mu_1\|u\|_p^p \nonumber\\
&-\frac{(q-2)N}{2q}\frac{(q-2)N-2}{2}\mu_2\|v\|_q^q\nonumber\\
&-\frac{(r_1+r_2-2)N}{2}\frac{(r_1+r_2-2)N-2}{2}\beta \int_{\R^N} |u|^{r_1} |v|^{r_2} \mathrm{d}x .
\end{align*}
On the other hand, for any $[u,v]\in \mathcal{P}$, we have
\begin{align*}
\|\nabla u\|_2^2+\|\nabla v\|_2^2=&\frac{(p-2)N}{2p}\mu_1\|u\|_p^p+\frac{(q-2)N}{2q}\mu_2\|v\|_q^q\nonumber\\
&+\frac{(r_1+r_2-2)N}{2}\beta \int_{\R^N} |u|^{r_1} |v|^{r_2} \mathrm{d}x.
\end{align*}
Hence,
\begin{align*}
\Psi''_{[u,v]}(1)=&-\left[(p-2)N-4\right] \frac{(p-2)N}{4p}\mu_1\|u\|_p^p-\left[(q-2)N-4\right]\frac{(q-2)N}{4q}\mu_2\|v\|_q^q\nonumber\\
&-\left[(r_1+r_2-2)N-4\right]\frac{(r_1+r_2-2)N}{4}\beta\int_{\R^N}|u|^{r_1}|v|^{r_2}\mathrm{d}x.
\end{align*}
So by $2+\frac{4}{N}<p,q,r_1+r_2<2^*$, and $[u,v]\neq [0,0]$, we see that $\Psi''_{[u,v]}(1)<0$.
\ep



\s{Rearrangement}\lab{sec:Rearrangement}
\renewcommand{\theequation}{5.\arabic{equation}}

For any $[u,v]\in \mathcal{P}_{(a,b)}$, let $u^*, v^*$ be the Schwartz symmetrization of $u$ and $v$. We refer to \cite{LiLo} for a definition and the main properties of Schwartz symmetrization. We have the following result.

\bl\lab{lemma:20210623-wl1}
Let $2+\frac{4}{N}<p,q,r_1+r_2<2^*$.
For any $[u,v]\in \mathcal{P}_{(a,b)}$, there exists a unique $t=t_{[u^*,v^*]}\in (0,1]$ such that $[t\star u^*, t\star v^*]\in \mathcal{P}_{(a,b)}$ and we have  $J[t\star u^*, t\star v^*]\leq J[u,v].$
\el
\bp
For any $[u,v]\in \mathcal{P}_{(a,b)}$, we have $[u,v]\neq [0,0]$. So by $\|u^*\|_2^2=\|u\|_2^2$ and $\|v^*\|_2^2=\|v\|_2^2$, we see that $[u^*,v^*]\in (D_a\times D_b)\backslash \{[0,0]\}$. Then by Corollary \ref{cro:uniqueness-t}, there exists a unique
$t=t_{[u^*,v^*]}>0$ such that $[t\star u^*, t\star v^*]\in \mathcal{P}_{(a,b)}$.
In particular, by the properties of rearrangement, we have
\begin{equation*}
\begin{cases}
\|\nabla u^*\|_2^2\leq \|\nabla u\|_2^2, \|\nabla v^*\|_2^2\leq \|\nabla v\|_2^2,\\
\|u^*\|_p^p=\|u\|_p^p, \|v^*\|_q^q=\|v\|_q^q,\\
\int_{\R^N} |u^*|^{r_1} |v^*|^{r_2}\mathrm{d}x\geq \int_{\R^N}|u|^{r_1}|v|^{r_2}\mathrm{d}x.
\end{cases}
\end{equation*}
This implies that $P[u^*,v^*]\leq P[u,v]=0$. Hence, by Corollary \ref{cro:uniqueness-t}, we obtain $t_{[u^*,v^*]}\leq 1$.
We remark that for any $[u,v]\neq [0,0]$, we have that $J[u^*, v^*]\leq J[u,v]$. Combined with the fact that
$$s\star u^*=(s\star u)^*, \quad \mbox{for all } s\in \R^+,$$
we finally obtain that
\begin{align*}
\max_{s>0}J[s\star u^*, s\star v^*]=J[t\star u^*, t\star v^*]=&J[(t\star u)^*, (t\star v)^*]
\leq J[t\star u, t\star v]\\
\leq&\max_{s>0}J[s\star u, s\star v]
=J[u,v].
\end{align*}
\ep

\begin{remark}\label{symmetric_ground_state}
Lemma \ref{lemma:20210623-wl1} implies, in particular, that if $C_{(a,b)}$ is achieved then it is achieved by a couple $[u,v] \in \mathcal{H}$ of Schwartz symmetric functions.
\end{remark}


\s{Existence of a special Palais-Smale sequence at the level $C_{(a,b)}$}\lab{sec:Palais_Smale}
\renewcommand{\theequation}{5.\arabic{equation}}

In this section we prove the existence of a particular Palais-Smale sequence for $J$ restricted to $D_a\times D_b$ at the level $C_{(a,b)}$.  \bigskip

Let us define the class of paths
\begin{equation*}
\G_{(a,b)} = \{g \in C([0,1], D_a\times D_b) : g(0) \in L_{(a,b)} \,\,  \mbox{and} \,\, P(g(1)) <0 \}
\end{equation*}
where $$L_{(a,b)} := \{ [u,v] \in \mathcal{H} : ||\nabla u ||_2^2 + ||\nabla u ||_2^2 \leq \frac{1}{2}\delta_{(a,b)} \}.$$

Note that $\G_{(a,b)} \neq \emptyset.$ Indeed $L_{(a,b)} \neq \emptyset$ and from Corollary \ref{cro:uniqueness-t} we known that the set where $P<0$ is non void.  \smallskip

Next we define
\begin{equation*}
\gamma_{(a,b)} = \inf_{g \in \G_{(a,b)}} \max_{t \in [0,1]} J(g(t)).
\end{equation*}
We claim that $$\gamma_{(a,b)} = C_{(a,b)}.$$ Indeed, combining Corollary \ref{cro:uniqueness-t} 	and Lemma \ref{lemma:tidu-xyj} we see that $P >0$ on $L_{(a,b)}$ and thus, by continuity, for any $g \in \G_{(a,b)}$ there exists a $t \in [0,1]$ such that $P(g(t)) =0$. It implies that $\gamma_{(a,b)} \geq C_{(a,b)}$. The reverse inequality $\gamma_{(a,b)} \leq C_{(a,b)}$ follows from Corollary \ref{cro:uniqueness-t}.  \smallskip

As a consequence of the equality $\gamma_{(a,b)} = C_{(a,b)}$ we deduce that
$$ \inf_{u \in \mathcal{P}_{(a,b)} } J(u) \geq \gamma_{(a,b)}.$$
At this point we observe that we are in the setting of \cite[Theorem 4.1]{Gh}. This theorem says to us that, for any sequence of paths $\{ g_n \} \subset \G_{(a,b)}$ such that,  as  $n \rightarrow \infty$,
\begin{equation}\label{mimpaths}
\max_{[0,1]} J(g_n(t)) \rightarrow C_{(a,b)},
\end{equation}
 there exists a sequence $\{[u_n,v_n] \} \subset \mathcal{H}$ such that, as $n \rightarrow \infty$,
\begin{itemize}
\item[(i)] $J[u_n, v_n] \rightarrow C_{(a,b)}$, \smallskip

\item[(ii)] $ J'[u_n,v_n]\Big|_{D_a \times D_b}\rightarrow 0$, \smallskip

\item[(iii)] $dist([u_n,v_n], \mathcal{P}_{(a,b)}) \rightarrow 0$,  \smallskip

\item[(iv)] $dist([u_n,v_n], g_n[0,1]) \rightarrow 0.$
\end{itemize}
\smallskip

Let $[u_n^*, v_n^*] \subset \mathcal{P}_{(a,b)}$ be a minimizing sequence for $C_{(a,b)}$. We know by Lemma \ref{lemma:20210623-wl1} that such a sequence exists. Now let the sequence of paths $\{g_n\} \subset \G_{(a,b)}$ be given by :
$$g_n(t) = [(s_1^n + t s_2^n) * u_n^*, (s_1^n + t s_2^n) * v_n^*]$$
where $s_1^n >0$ is sufficiently small and $s_2^n >0$ sufficiently large in order to insure that $g_n \in \G_{(a,b)}$. We know from Corollary \ref{cro:def-Cab} that the sequence $\{g_n\}$ satisfy \eqref{mimpaths}.
 \smallskip

From (i)-(ii) we deduce that $\{[u_n,v_n] \} \subset \mathcal{H}$ is a Palais-Smale sequence for $J$ restricted to $D_a \times D_b$ at the level $C_{(a,b)}$. The property (iii), combined with the fact that $J$ is coercive on $\mathcal{P}_{(a,b)}$, see
Corollary \ref{cro:20210623-c1}, and that $J'$ and $P'$ take bounded sets into bounded set implies that $\{[u_n,v_n] \} \subset \mathcal{H}$ is bounded and that $P[u_n,v_n] \rightarrow 0$.  Now, by compactness, for each $n \in \N$, there exists $s_n \in [s_1^n, s_1^n+s_2^n]$ such that
\begin{equation}\label{distance_to_paths}
dist \Big([u_n, v_n], g_n([0,1])\Big) = \Big|\Big| [u_n,v_n] - [s_n * u_n^*, s_n * v_n^*] \Big|\Big|.
\end{equation}
Up to a subsequence we can assume that $[s_n * u_n^*, s_n * v_n^*] \rightharpoonup u $ weakly in $\mathcal{H}$ and
$[s_n * u_n^*, s_n * v_n^*] \rightarrow u $
 strongly in $L^{\eta}(\R^N) \times L^{\eta}(\R^N)$ for all $2 < \eta < 2^*$. For $N \geq 2$ this strong convergence is a consequence of the fact that the functions are radial. When $N=1$ it follows from the fact that the functions are nonincreasing  with respect to the origin, see  \cite[Proposition 1.7.1]{Cazenave2003}.

Now, in view of \eqref{distance_to_paths}, these convergence properties hold for the Palais-Smale sequence $\{[u_n,v_n] \} \subset \mathcal{H}.$

Summarizing we have obtained a bounded sequence $\{[u_n,v_n] \} \subset D_a\times D_b$ such that
\begin{itemize}
\item[(i)] $J[u_n,v_n]\rightarrow C_{(a,b)}$,
\item[(ii)] $ J'[u_n,v_n] + \lambda_{1,n}[u_n,0] + \lambda_{2,n}[0,v_n] \to 0, \quad  \mbox{in } \mathcal{H}^{-1}$ for some real sequences $\{\lambda_{1,n}\}$ and $\{\lambda_{2,n}\}$,
\item[(iii)] $ P[u_n, v_n] \rightarrow 0$.
\item[(iv)] $[u_n,v_n] \rightharpoonup [u,v] $   weakly in  $\mathcal{H}$ and  $[u_n,v_n] \rightarrow [u,v] $ strongly in  $L^{\eta}(\R^N) \times L^{\eta}(\R^N)$ for all $2 < \eta < 2^*$. In addition $u \geq 0$ and $v\geq 0$.
\end{itemize}

\s{Estimation of $C_{(a,b)}$}\label{sec:estimation}
\renewcommand{\theequation}{6.\arabic{equation}}
In view of the definitions of $z_{p,\mu,a}$ by \eqref{def:zpua} and $\beta_{p,\mu,a,N,r}$ by \eqref{def:-best-const}, a direct calculation shows that
\begin{equation*}
\beta_{p,\mu,a,N,r}=\frac{1}{2}\lambda_{p,\mu,a}^{\frac{p-2-r}{p-2}} \mu^{\frac{r}{p-2}}
\inf_{h\in H^1(\R^N)\backslash\{0\}}\frac{\int_{\R^N}|\nabla h|^2 \mathrm{d}x}{\int_{\R^N}U_p^r h^2 \mathrm{d}x}.
\end{equation*}
Substituting the value of $\lambda_{p,\mu,a}$ (see \eqref{eq:20220905-ze1}), we obtain that
\begin{equation*}
\lab{eq:20221212-be1}
\beta_{p,\mu,a,N,r}=\frac{1}{2}\mu^{\frac{rN-4}{N(p-2)-4}} \,
\|U_p\|_{2}^{\frac{4(p-2-r)}{N(p-2)-4}} a^{-\frac{2(p-2-r)}{N(p-2)-4}} \inf_{h\in H^1(\R^N)\backslash\{0\}}\frac{\int_{\R^N}|\nabla h|^2 \mathrm{d}x}{\int_{\R^N}U_p^r h^2 \mathrm{d}x}.
\end{equation*}

\bl\lab{lemma:20220907-l1}
For $N=1,2$, $p\in (2,+\infty), \mu>0, a>0, r>0$, it holds that
$\beta_{p,\mu,a,N,r}=0$.
\el
\bp
For the case of $N=1$, let $h(x)\in H^1(\R)$ such that $h(x)\equiv 1$ for $|x|\leq 1$, $h(x)\equiv 0$ for $|x|\geq 2$ and $0\leq h(x)\leq 1$.  Define
$\displaystyle h_n(x):=h\left(\frac{x}{n}\right)$.
Then $\|\nabla h_n\|_2^2=\frac{1}{n}\|\nabla h\|_2^2\rightarrow 0$ as $n\rightarrow \infty$.
While
$$\int_{\R}U_p^r(x) h_n^2(x) \mathrm{d}x \geq \int_{|x|\leq n} U_p^r (x)\mathrm{d}x \rightarrow \int_{\R} U_p^r(x)\mathrm{d}x>0.$$
So
$$\inf_{0\neq h\in H^1(\R)}\frac{\|\nabla h\|_2^2}{\int_{\R}U_p^r(x) h(x)^2 \mathrm{d}x}=0$$
and thus $\beta_{p,\mu,a,N,r}=0$ for $N=1$.

For $N=2$, we remark that the embedding $H^1(\R^2)\hookrightarrow L^\infty(\R^2)$ fails, so there exists a sequence $\{\phi_n\}\subset H^1(\R^N)$ with $\|\phi_n\|_{H^1}=1$, but $\|\phi_n\|_\infty\rightarrow \infty$ as $n\rightarrow \infty$. In particular, by the property of rearrangement, without loss of generality, we can assume that $\phi_n=\phi_n^*$. So we have that
$$\phi_n(x)=\phi_n(|x|)~\hbox{and}~\phi_n(0)=\|\phi_n\|_\infty.$$
Now, we set
$$h_n(x)=\begin{cases}
\phi_n(0)\quad &\hbox{if}~|x|\leq n,\\
\phi_n(|x|-n)\quad &\hbox{if}~|x|\geq n.
\end{cases}$$
Then $\|\nabla h_n\|_2^2=\|\nabla \phi_n\|_2^2\leq \|\phi_n\|_{H^1}^{2}=1$ and
\begin{align*}
&\int_{\R^2}U_p^r(x)h_n^2(x) \mathrm{d}x\geq \int_{|x|\leq n} \phi_n^2(0) U_p^r(x) \mathrm{d}x
=\|\phi_n\|_\infty^2 \int_{|x|\leq n} U_p^r(x)\mathrm{d}x
\rightarrow \infty~\hbox{as}~n\rightarrow \infty.
\end{align*}
So we also have
$$\inf_{0\neq h\in H^1(\R^2)}\frac{\|\nabla h\|_2^2}{\int_{\R^2}U_p^r(x) h^2(x)} \mathrm{d}x=0,$$
and thus $\beta_{p,\mu,a,N,r}=0$ also holds for $N=2$.
\ep

\bl\lab{20220905-l1}
For $N\geq 3$, it holds that $\beta_{p,\mu,a,N,r}>0$. In particular, for $r>\frac{4}{N}$ and fixed $2+\frac{4}{N}<p<2^*,a>0$, we have
\beq\lab{eq:20220905-xe1}
\lim_{\mu\rightarrow \infty}\beta_{p,\mu,a,N,r}=+\infty.
\eeq
Furthermore, for $2+\frac{4}{N}<p<2^*, \mu>0, r>1$ fixed, we have that
\begin{itemize}
\item[(i)] if $p<r+2<2^*$, then $\displaystyle \lim_{a\rightarrow 0}\beta_{p,\mu,a,N,r}=0$ and $\displaystyle \lim_{a\rightarrow \infty}\beta_{p,\mu,a,N,r}=+\infty$.
\item[(ii)] if $r+2=p$, then there exists $\tilde{C}_{p,\mu}<\tilde{D}_{p,\mu}$ independent of $a$, such that
$$\tilde{C}_{p,\mu}\leq \beta_{p,\mu,a,N,r} \leq \tilde{D}_{p,\mu}, \forall a>0.$$
\item[(iii)] if $2+r<p$, then $\displaystyle \lim_{a\rightarrow 0}\beta_{p,\mu,a,N,r}=+\infty$ and $\displaystyle \lim_{a\rightarrow \infty}\beta_{p,\mu,a,N,r}=0$.
\end{itemize}
\el
\bp
Noting that $U_p\in L^\infty(\R^N)$ and decays exponentially to $0$ as $|x|\rightarrow +\infty$,
by H\"older inequality, we have that
\begin{equation*}
\int_{\R^N}U_p^{r} |h|^2 \mathrm{d}x\leq \|U_p\|_{\frac{2^* r}{2^*-2}}^{r}  \|h\|_{2^*}^{2}.
\end{equation*}
So by the critical Sobolev inequality, we obtain that
\begin{align}\lab{eq:20220905-zbe1}
\beta_{p,\mu,a,N,r}\geq & \frac{1}{2}\mu^{\frac{rN-4}{N(p-2)-4}} \,
\|U_p\|_{2}^{\frac{4(p-2-r)}{N(p-2)-4}} a^{-\frac{2(p-2-r)}{N(p-2)-4}}\frac{\|\nabla h\|_2^2}{\|U_p\|_{\frac{2^* r}{2^*-2}}^{r}  \|h\|_{2^*}^{2}}\nonumber\\
\geq& \frac{1}{2}\mu^{\frac{rN-4}{N(p-2)-4}} \,
\|U_p\|_{2}^{\frac{4(p-2-r)}{N(p-2)-4}} a^{-\frac{2(p-2-r)}{N(p-2)-4}}\|U_p\|_{\frac{2^* r}{2^*-2}}^{-r} S_N,
\end{align}
where $S_N >0$ is the sharp constant in the critical Sobolev inequality.
So it is easy to see that \eqref{eq:20220905-xe1} holds.

On the other hand, by taking $h=U_p$, we obtain that
\beq\lab{eq:20220905-zbe2}
\beta_{p,\mu,a,N,r}\leq
\frac{1}{2}\mu^{\frac{rN-4}{N(p-2)-4}}
\|U_p\|_{2}^{\frac{4(p-2-r)}{N(p-2)-4}} a^{-\frac{2(p-2-r)}{N(p-2)-4}} \frac{\|\nabla U_p\|_2^2}{\|U_p\|_{r+2}^{r+2}} .
\eeq
Then by \eqref{eq:20220905-zbe1} and \eqref{eq:20220905-zbe2}, it is easy to see that the conclusions (i)-(iii) above hold.
\ep

\bl\lab{lemma:20221001-wl1}
\begin{itemize}
\item[(i)]if $1<r_1<2$, then $C_{(a,b)}<m_{q,\mu_2,b}$. If $r_1=2$, then  $C_{(a,b)}<m_{q,\mu_2,b}$ provided $\beta>\beta_{q,\mu_2,b,N,r_2}$.
\item[(ii)]if $1<r_2<2$, then $C_{(a,b)}<m_{p,\mu_1,a}$. If $r_2=2$, then $C_{(a,b)}<m_{p,\mu_1,a}$ provided $\beta>\beta_{p,\mu_1,a,N,r_1}$.
\end{itemize}
\el
\bp
We only prove (ii), the proof of (i) being directly similar. We write $z:=z_{p,\mu_1,a}$ for simplicity. For any $h\in H^1(\R^N)$ with $\|h\|_2^2=1$, we have that $[z, sh]\in D_a\times D_b$  provided $|s|<\sqrt{b}$.
For any $s$, there exists a unique $t=t(s)>0$ such that
$[t\star z, t\star sh]\in \mathcal{P}_{(a,b)}$, where $t=t(s)$ is determined by
\begin{align*}
\|\nabla z\|_2^2+s^2 \|\nabla h\|_2^2=&\frac{(p-2)N}{2p}\mu_1\|z\|_p^p \, t^{\frac{(p-2)N-4}{2}}\nonumber\\
&+\frac{(q-2)N}{2q}\mu_2\|h\|_q^q \,  s^q \, t^{\frac{(q-2)N-4}{2}}\nonumber\\
&+\frac{(r_1+r_2-2)N}{2}\beta \Big(\int_{\R^N}|z|^{r_1}|h|^{r_2}\mathrm{d}x\Big) s^{r_2} \,  t^{\frac{(r_1+r_2-2)N-4}{2}}.
\end{align*}
From the Implicit Function Theorem we know that $t(s)\in C^1$ locally around $s=0$.
 Taking the derivative with respect to $s$ at both sides and we obtain that
 \begin{align*}
 t'(s)=\frac{P_h(s)}{Q_h(s)},
 \end{align*}
 where
 \begin{align*}
 P_h(s):=&2\|\nabla h\|_2^2 s -\frac{(q-2)N}{2}\mu_2 \, \|h\|_q^q \, |s|^{q-2}s \, t^{\frac{(q-2)N-4}{2}}\nonumber\\
 &-\frac{(r_1+r_2-2)N}{2}\beta \, r_2 \Big(\int_{\R^N}|z|^{r_1} |h|^{r_2}\mathrm{d}x \Big) |s|^{r_2-2}s \, t^{\frac{(r_1+r_2-2)N-4}{2}}
 \end{align*}
 and
 \begin{align*}
 Q_h(s):=&\frac{(p-2)N-4}{2}\frac{(p-2)N}{2p}\mu_1 \, \|z\|_p^p \,  t^{\frac{(p-2)N-6}{2}}\nonumber\\
 &+\frac{(q-2)N-4}{2}\frac{(q-2)N}{2q}\mu_2 \, \|h\|_q^q \,  s^q \, t^{\frac{(q-2)N-6}{2}}\nonumber\\
 &+\frac{(r_1+r_2-2)N-4}{2}\frac{(r_1+r_2-2)N}{2}\beta \Big(\int_{\R^N}|z|^{r_1}|h|^{r_2}\mathrm{d}x \Big) s^{r_2} \, t^{\frac{(r_1+r_2-2)N-6}{2}}.
 \end{align*}

 {\bf Case of $1<r_2<2$:} In such a case, for $|s|$ small,
 $$P_h(s)=-\frac{(r_1+r_2-2)N}{2}\beta r_2 \Big(\int_{\R^N}|z|^{r_1} |h|^{r_2}\mathrm{d}x \Big) |s|^{r_2-2}s (1+o(1))$$
 and
 $$Q_h(s)=\frac{(p-2)N-4}{2} \, \frac{(p-2)N}{2p}\mu_1\|z\|_p^p (1+o(1)).$$
 So
 \begin{equation*}
 t'(s)=-M_h r_2 |s|^{r_2-2}s (1+o(1)),
 \end{equation*}
 where
 \begin{equation*}
 M_h:=\frac{2(r_1+r_2-2)N}{(p-2)N-4} \, \frac{\beta \int_{\R^N}|z|^{r_1} |h|^{r_2}\mathrm{d}x}{\|\nabla z\|_2^2}.
 \end{equation*}
 Then
 \begin{equation*}
 t(s)=1-M_h|s|^{r_2} (1+o(1))
 \end{equation*}
 and for any $\tau>0$,
 \begin{equation*}
 t(s)^\tau=1-\tau M_h |s|^{r_2} (1+o(1)).
 \end{equation*}
 Noting that for any $[\phi,\psi]\in \mathcal{P}$, it holds that
\begin{align*}
J[\phi,\psi]=&\frac{(p-2)N-4}{4} \,\frac{\mu_1}{p}\|\phi\|_p^p +\frac{(q-2)N-4}{4} \,\frac{\mu_2}{q}\|\psi\|_q^q \\
&+\frac{(r_1+r_2-2)N-4}{4}\beta \int_{\R^N}|\phi|^{r_1}|\psi|^{r_2}\mathrm{d}x.
\end{align*}
Then we have
 \begin{align*}
 J[t(s)\star z, t(s)\star (sh)]-J[z,0]
 &=\frac{1}{2}\|\nabla z\|_2^2 \, (t^2(s)-1) -\frac{\mu_1}{p} \, \|z\|_p^p \, (t(s)^{\frac{(p-2)N}{2}}-1)\\
 &+\frac{1}{2}\|\nabla h\|_2^2 \,  s^2 \, t^2(s) -\frac{\mu_2}{q} \, \, \|h\|_q^q \, s^q \, t(s)^{\frac{(q-2)N}{2}}\\
 &-\beta \Big(\int_{\R^N}|z|^{r_1}|h|^{r_2}\mathrm{d}x\Big) |s|^{r_2} \, t(s)^{\frac{(r_1+r_2-2)N}{2}}\\
 =&-\|\nabla z\|_2^2 \, M_h \, |s|^{r_2} o(1) -\beta \Big(\int_{\R^N}|z|^{r_1}|h|^{r_2}\mathrm{d}x \Big) |s|^{r_2}(1+o(1))\\
 =&-\beta \Big(\int_{\R^N}|z|^{r_1}|h|^{r_2}\mathrm{d}x \Big) |s|^{r_2}(1+o(1)).
 \end{align*}
Hence, if $1<r_2<2$, by taking $s$ close to $0$, we obtain that
\begin{equation*}
C_{(a,b)}\leq J[t(s)\star z, t(s)\star (sh)]<J[z,0]=m_{p,\mu_1,a}.
\end{equation*}

{\bf Case of $r_2=2$}: In such a case, we have, for $s>0$ small,
 $$P_h(s)=\left[2\|\nabla h\|_2^2-N\beta r_1 \int_{\R^N}|z|^{r_1} |h|^{2}\mathrm{d}x\right] s (1+o(1))$$
 and
\begin{equation*}
 t'(s)=\bar{M}_h s(1+o(1))
\end{equation*}
where
\begin{equation*}
\bar{M}_h:=\frac{2\|\nabla h\|_2^2-N\beta r_1 \int_{\R^N}|z|^{r_1} |h|^{2}\mathrm{d}x}{\frac{(p-2)N-4}{2}\|\nabla z\|_2^2}.
\end{equation*}
Then
$$t(s)=1+\frac{1}{2}\bar{M}_h \,s^2(1+o(1))$$
and for any $\tau>0$,
$$t(s)^\tau=1+\frac{\tau}{2}\bar{M}_h \,s^2(1+o(1)).$$
So by a similar argument as in the  Case of $1<r_2<2$, we obtain that
\begin{align*}
J[t(s)\star z, t(s)\star (sh)]-J[z,0]
 =&\frac{1}{2}\|\nabla z\|_2^2 \, (t(s)^2-1) -\frac{1}{p}\, \mu_1\|z\|_p^p (t(s)^{\frac{(p-2)N}{2}}-1)\\
 &+\frac{1}{2}\|\nabla h\|_2^2 \, s^2 \, t(s)^2 -\frac{1}{q}\mu_2\,  \|h\|_q^q \, s^q \, t(s)^{\frac{(q-2)N}{2}}\\
 &-\beta \Big( \int_{\R^N}|z|^{r_1}|h|^{2}\mathrm{d}x \Big) |s|^{2} t(s)^{\frac{r_1N}{2}}\\
 =&\left[\frac{1}{2}\|\nabla h\|_2^2 -\beta  \int_{\R^N}|z|^{r_1}|h|^{2}\mathrm{d}x\right]s^2(1+o(1)).
\end{align*}
Then for any $\beta>\beta_{p,\mu_1,a,N,r_1}$, there exists some $h\in H^1(\R^N)$ such that
$$\frac{1}{2}\|\nabla h\|_2^2 -\beta \int_{\R^N}|z|^{r_1}|h|^{2}\mathrm{d}x<0,$$
and for such a $h$ and $|s|$ small, we have
$J[t(s)\star z, t(s)\star (sh)]-J[z,0]<0$. Hence,
$$C_{(a,b)}<m_{p,\mu_1,a}.$$
\ep

\s{Convergence of our Palais-Smale sequence and Proof of Theorem \ref{th:main-t1}}\label{sec:PS-sequence}
\renewcommand{\theequation}{7.\arabic{equation}}

Theorem \ref{th:main-t1} will be a direct consequence of the following result.

\bt\lab{thm:compactness}
Let $1\leq N\leq 4$, $2+\frac{4}{N}<p,q,r_1+r_2<2^*$ and $\mu_1,\mu_2,\beta, a,b\in \R^+$.  Let
 $\{[u_n,v_n] \} \subset D_a\times D_b$ be a bounded sequence such that
\begin{itemize}
\item[(i)] $J[u_n,v_n]\rightarrow C_{(a,b)}$,
\item[(ii)] $ J'[u_n,v_n] + \lambda_{1,n}[u_n,0] + \lambda_{2,n}[0,v_n] \to 0, \quad  \mbox{in } \mathcal{H}^{-1}$ for some real sequences $\{\lambda_{1,n}\}$ and $\{\lambda_{2,n}\}$,
\item[(iii)] $ P[u_n, v_n] \rightarrow 0$,
\item[(iv)] $[u_n,v_n] \rightharpoonup [u,v] $   weakly in  $\mathcal{H}$ and  $[u_n,v_n] \rightarrow [u,v] $ strongly in  $L^{\eta}(\R^N) \times L^{\eta}(\R^N)$ for all $2 < \eta < 2^*$. In addition $u \geq 0$ and $v\geq 0$.
\end{itemize}
Assume that
\begin{equation}\label{strict_inequality}
C_{(a,b)}<\min \{m_{p,\mu_1,a}, m_{q,\mu_2,b}\}.
\end{equation}
Then, up to a subsequence, $[u_n,v_n]\rightarrow [u,v]$ in $\mathcal{H}$ and  $[u,v]\in S_a\times S_b$.
\et

\bp
Without loss of generality, we may assume that $[u_n,v_n]\neq [0,0]$ for all $n\in \N$ since $J[u_n,v_n]\rightarrow C_{a,b}>0=J[0,0]$, see Corollary \ref{cro:def-Cab}.
First note that (ii) can be rewrite as : there exists sequences $\{\lambda_{1,n}\}$ and $\{\lambda_{2,n}\}$ such that
\begin{equation*}
\begin{cases}
-\Delta u_n+\lambda_{1,n}u_n-\mu_1 u_{n}^{p-1}-\beta r_1 u_{n}^{r_1-1}v_{n}^{r_2}=o(1),\\
-\Delta v_n+\lambda_{2,n}v_n-\mu_2 v_{n}^{q-1}-\beta r_2 u_{n}^{r_1}v_{n}^{r_2-1}=o(1).
\end{cases}
\end{equation*}
From the boundedness of $\{[u_n,v_n]\} \subset \mathcal{H}$,  it follows that $\{\lambda_{1,n}\}$ and $\{\lambda_{2,n}\}$ are also bounded and we can assume that $\lambda_{1,n}\rightarrow \lambda_1$ and $\lambda_{2,n}\rightarrow \lambda_2$ for some $\lambda_1, \lambda_2$.

We claim that
\beq\lab{eq:20220830-we1}
u\neq 0, v\neq 0.
\eeq
If not, without loss of generality, we may assume that $u=0$. Then we have that $v\neq 0$. If not, from $P[u_n,v_n] \to 0$ we get that
\begin{align*}
\|\nabla u_n\|_2^2+\|\nabla v_n\|_2^2=&\frac{(p-2)N}{2p}\mu_1\|u_n\|_p^p+\frac{(q-2)N}{2q}\mu_2\|v_n\|_q^q\nonumber\\
&+\frac{(r_1+r_2-2)N}{2}\beta \int_{\R^N} |u_n|^{r_1} |v_n|^{r_2} \mathrm{d}x + o(1),\nonumber\\
=&o(1),
\end{align*}
where we have used that
$$
\|u_n\|_p^p=o(1), \|v_n\|_q^q=o(1) \mbox{ and } \int_{\R^N}u_{n}^{r_1}v_{n}^{r_2}\mathrm{d}x=o(1),
$$
by the H\"older inequality and Gagliardo-Nirenberg inequality. This contradicts Lemma \ref{lemma:tidu-xyj}.
We deduce that
$$0<\|v\|_2^2\leq \liminf_{n\rightarrow \infty}\|v_n\|_2^2.$$
Thus one can see that $v$ is a non-trivial (non-negative) solution to
\beq\lab{eq:20210623-wze1}
-\Delta v+\lambda_2 v=\mu_2 v^{q-1}\;\hbox{in}~\R^N, \, v\in H^1(\R^N).
\eeq
Recalling the Pohozaev identity
\begin{equation*}
\|\nabla v\|_2^2=\frac{(q-2)N}{2q}\mu_2\|v\|_q^q,
\end{equation*}
we deduce, from \eqref{eq:20210623-wze1}, that
$$\lambda_2\|v\|_2^2=\frac{2q-(q-2)q}{2q}\mu_2 \|v\|_q^q>0.$$
Hence, $\lambda_2>0$. So by
\begin{align*}
\|\nabla v_n\|_2^2+\lambda_{2,n}\|v_n\|_2^2=&\mu_2\|v_n\|_q^q+\beta r_2\int_{\R^N}u_{n}^{r_1} v_{n}^{r_2}\mathrm{d}x+o(1)\\
=&\mu_2\|v\|_q^q+o(1)\\
=&\|\nabla v\|_2^2+\lambda_{2}\|v\|_2^2+o(1),
\end{align*}
we have that $v_n\rightarrow v$ in $H^1(\R^N)$.
Now using that $P[u_n,v_n] \to 0,$ and $ u_n\rightarrow 0$ in $L^\eta(\R^N)$ for all $2<\eta<2^*$, it is easy to prove that
$u_n\rightarrow 0$ in $D_{0}^{1,2}(\R^N)$.
So
$$ J[0,v]=\lim_{n\rightarrow \infty}J[u_n,v_n]=C_{(a,b)}.$$
Put $\delta:=\|v\|_2^2$, then $\delta\in (0,b]$.
By Lemma \ref{lemma:20210622-l1}, one has that
$v=z_{q,\mu_2,\delta}$ and $C_{(a,b)}=J[0,v]=m_{q,\mu_2,\delta}$.
However, by Lemma \ref{lemma:20210622-l2}-(i), and \eqref{strict_inequality} we have
$$m_{q,\mu_2,\delta}\geq  m_{q,\mu_2,b}> C_{(a,b)},$$
 a contradiction. This ends the proof of \eqref{eq:20220830-we1} and thus
\begin{equation*}
0<\|u\|_2^2\leq \liminf_{n\rightarrow \infty}\|u_n\|_2^2 \quad \mbox{and} \quad 0<\|v\|_2^2\leq \liminf_{n\rightarrow \infty}\|v_n\|_2^2.
\end{equation*}
Now, from the convergence properties (iv), we deduce that
 $[u,v]$ is a non-trivial, solution to \eqref{eq:main-equation}. In particular it holds that $P[u,v]=0$ and  by
\begin{align*}
&\|\nabla u\|_2^2+\|\nabla v\|_2^2\leq o(1)+\|\nabla u_n\|_2^2+\|\nabla v_n\|_2^2\\
=&o(1)+\frac{(p-2)N}{2p}\mu_1\|u_n\|_p^p+\frac{(q-2)N}{2q}\mu_2\|v_n\|_q^q+\frac{(r_1+r_2-2)N}{2}\beta \int_{\R^N} |u_n|^{r_1} |v_n|^{r_2} \mathrm{d}x\\
=&o(1)+\frac{(p-2)N}{2p}\mu_1\|u\|_p^p+\frac{(q-2)N}{2q}\mu_2\|v\|_q^q+\frac{(r_1+r_2-2)N}{2}\beta \int_{\R^N} |u|^{r_1} |v|^{r_2} \mathrm{d}x\\
=&o(1)+\|\nabla u\|_2^2+\|\nabla v\|_2^2,
\end{align*}
we obtain that $u_n\rightarrow u, v_n\rightarrow v$ in $D_{0}^{1,2}(\R^N)$. Hence, $\displaystyle J[u,v]=\lim_{n\rightarrow \infty}J[u_n,v_n]=C_{(a,b)}.$

\noindent
Recalling \eqref{eq:main-equation} and $1\leq N\leq 4$, we conclude, using \cite[Lemma A.2]{ikoma2014compactness}, that $\lambda_1>0,\lambda_2>0$.
Then by
\begin{align*}
\|\nabla u_n\|_2^2+\lambda_{1,n}\|u_n\|_2^2=&\mu_1\|u_n\|_p^p+\beta r_1\int_{\R^N}u_{n}^{r_1} v_{n}^{r_2}\mathrm{d}x+o(1)\\
=&\mu_1\|u\|_p^p+\beta r_1\int_{\R^N}u^{r_1} v^{r_2}\mathrm{d}x+o(1)\\
=&\|\nabla u\|_2^2+\lambda_{1}\|u\|_2^2+o(1),
\end{align*}
we obtain that $u_n\rightarrow u$ in $H^1(\R^N)$. Similarly, we can prove that $v_n\rightarrow v$ in $H^1(\R^N)$.
To obtain that $[u,v]\in S_a\times S_b$ and completed the proof we shall use Lemma \ref{lemma:20220831-l1} below.
\ep

\bl\lab{lemma:20220831-l1}
Let $[u,v], \lambda_1,\lambda_2$ be given by Theorem \ref{thm:compactness}. Since $\lambda_1>0$ we have that $u\in S_a$. Similarly,  $\lambda_2>0$ imply that $v\in S_b$.
\el
\bp
Suppose that $\lambda_1>0$, we shall prove that $u\in S_a$. If not, $\delta:=\|u\|_2^2\in (0, a)$. Then for $s>0$ small enough, we still have that $[(1+s)u, v]\in (D_a\times D_b)\backslash\{[0,0]\}$. For any given $s$, by Corollary \ref{cro:uniqueness-t}, there exists a unique $t=t(s)>0$ such that
$[t\star (1+s)u, t\star v]\in \mathcal{P}_{(a,b)}$.
Precisely, $t=t(s)$ is determined by
{\allowdisplaybreaks
\begin{align}
(1+s)^2\|\nabla u\|_2^2+\|\nabla v\|_2^2=&\frac{(p-2)N}{2p}\mu_1\|u\|_p^p \, t^{\frac{(p-2)N-4}{2}}(1+s)^p\nonumber\\
&+\frac{(q-2)N}{2q}\mu_2\|v\|_q^q \, t^{\frac{(q-2)N-4}{2}}\nonumber\\
&+\frac{(r_1+r_2-2)N}{2}\beta \Big(\int_{\R^N} |u|^{r_1} |v|^{r_2} \mathrm{d}x \Big) t^{\frac{(r_1+r_2-2)N-4}{2}}(1+s)^{r_1}.\nonumber
\end{align}}
From the Implicit Function Theorem we have that $t(s)\in C^1$. Then, since
\begin{align}
J[t\star (1+s)u, t\star v]=&\frac{1}{2}[\|\nabla u\|_2^2(1+s)^2+\|\nabla v\|_2^2]t^2-\frac{\mu_1}{p}\|u\|_p^p (1+s)^p \, t^{\frac{(p-2)N}{2}}\nonumber\\
&-\frac{\mu_2}{q}\|v\|_q^q \, t^{\frac{(q-2)N}{2}}-\beta \Big(\int_{\R^N} |u|^{r_1} |v|^{r_2} \mathrm{d}x \Big) t^{\frac{(r_1+r_2-2)N}{2}}(1+s)^{r_1},\nonumber
\end{align}
we have
{\allowdisplaybreaks
\begin{align*}
\frac{d}{ds}J[t\star (1+s)u, t\star v]
=&\|\nabla u\|_2^2 (1+s) t^2-\mu_1\|u\|_p^p (1+s)^{p-1} \,  t^{\frac{(p-2)N}{2}}\nonumber\\
&-\beta r_1\int_{\R^N}|u|^{r_1}|v|^{r_2}\mathrm{d}x (1+s)^{r_1-1} \,  t^{\frac{(r_1+r_2-2)N}{2}}\nonumber\\
&+\left(\|\nabla u\|_2^2 (1+s)^2+\|\nabla v\|_2^2\right)t \, t'-\frac{(p-2)N}{2p}\mu_1\|u\|_p^p (1+s)^p \,  t^{\frac{(p-2)N-2}{2}} \,t'\nonumber\\
&-\frac{(q-2)N}{2q}\mu_2\|v\|_q^q t^{\frac{(q-2)N-2}{2}} \, t'\nonumber\\
&-\frac{(r_1+r_2-2)N}{2}\beta \Big(\int_{\R^N}|u|^{r_1}|v|^{r_2}\mathrm{d}x \Big) (1+s)^{r_1}\, t^{\frac{(r_1+r_2-2)N-2}{2}} \, t',\nonumber
\end{align*}
}
here $t=t(s), t'=t'(s)$.
Putting $s=0$ and noting that $t(0)=1$, we have that
\begin{align*}
\frac{d}{ds}J[t(s)\star (1+s)u, t(s)\star v]\Big|_{s=0}
=&\big[\|\nabla u\|_2^2-\mu_1\|u\|_p^p-\beta r_1\int_{\R^N}|u|^{r_1}|v|^{r_2}\mathrm{d}x \big]+P[u,v]t'(0)\\
=&-\lambda_1\|u\|_2^2.
\end{align*}
Then
$$C_{(a,b)}\leq J[t(s)\star (1+s)u, t(s)\star v]<J[u,v]=C_{(a,b)}, \hbox{ for $s>0$ small enough},$$
a contradiction.
Similarly, using that $\lambda_2>0$, we can prove that $v\in S_b$.
\ep

At this point we can give the \smallskip

{\bf Proof of Theorem \ref{th:main-t1}:}
Under the assumptions, by Lemma \ref{lemma:20221001-wl1}, we have that
$$ C_{(a,b)}<\min \{m_{p,\mu_1,a}, m_{q,\mu_2,b}\}.$$ Let $\{[u_n,v_n]\} \subset \mathcal{H}$ be a Palais-Smale sequence as obtained in Section \ref{sec:Palais_Smale}. By Theorem \ref{thm:compactness}, $[u_n, v_n] \rightarrow [u,v] \in S_a \times S_b$ and, in particular, $[u,v] \in \mathcal{P}_{a,b}$ and $J[u,v] = C_{(a,b)}$. Invoking Remark \ref{symmetric_ground_state} we deduce that $J$ restricted to $\mathcal{P}_{(a,b)}$ admits a minimum, say $[\overline{u}, \overline{v}]$ which consists in
Schwarz symmetric functions. The fact that this minimum is a critical point of $J$ constrained to $D_a \times D_b$ follows from Lemmas \ref{lemma:20210623-l1} and \ref{lemma:20210623-l2}. As in the proof of Theorem \ref{thm:compactness} one deduces that the associated Lagrange multipliers $\overline{\lambda}_1$ and $\overline{\lambda}_2$  are both strictly positive. At this point we have proved that
$(\overline{\lambda}_1,\overline{\lambda}_2,\overline{u},\overline{v})\in \R^2\times \mathcal{H}$ is a ground state to Problem \eqref{eq:20220902-maine1} having the required symmetry properties.
\hfill$\Box$


\end{document}